\input amssym.tex

\mag=1200 
\hsize=130mm  \vsize=180mm  \voffset 5mm
\lineskiplimit=10pt \lineskip=10pt
\tolerance=10000 \pretolerance=1000 \parindent=0mm \raggedright

 \font\uc =cmr10 at  15pt

 \font\dz =cmr10 at  20pt 
 \font\dc =cmr10 at  25pt 

\font\tengo=eufm10

\font\sevengo=eufm7

\font\fivego=eufm5

\font\tenbb=msbm7 at 10pt

\font\sevenbb=msbm7   
 
\font\fivebb=msbm5


\newfam\gofam  \textfont\gofam=\tengo
\scriptfont\gofam=\sevengo   \scriptscriptfont\gofam=\fivego
\def\go{\fam\gofam\tengo}

\newfam\bbfam  \textfont\bbfam=\tenbb
\scriptfont\bbfam=\sevenbb   \scriptscriptfont\bbfam=\fivebb
\def\bb{\fam\bbfam\tenbb}


\def\ind{\hskip 1em\relax}

\def\Im{\mathop{\rm Im}\nolimits}

\def\Ker{\mathop{\rm Ker}\nolimits}

\def\End{\mathop{\rm End}\nolimits}

\def\dim{\mathop{\rm dim}\nolimits}
\def\deg{\mathop{\rm deg}\nolimits}

\let\f=\varphi

\let\mbox=\hbox
\let\then=\Longrightarrow

\let\dsp=\displaystyle

\def\Ext{\mathop{\rm Ext}\nolimits}
\def\R{{\bb R}} \def\C{{\bb C}}  \def\Z{{\bb Z}}

\def\cl{\centerline}

\def\{{\lbrace}  
\def\}{\rbrace}   
\def\({\langle}  
\def\){\rangle}
\def\¾{\leq} 
\def\"{\geq}
\def\og{\leavevmode\raise.3ex\hbox{$\scriptscriptstyle\langle\!\langle$}}
\def\fg{\leavevmode\raise.3ex\hbox{$\scriptscriptstyle\,\rangle\!\rangle$}}
\def\[{\lbrack} \def\]{\rbrack}
\def\l{\ell}
\def\arrow{\rightarrow}
\def\iso{\buildrel\sim\over{\arrow}}
\def\.{\bullet}
\def\bs{\bigskip}

\def\mono{\rightarrowtail}
\def\epi{\twoheadrightarrow}
\def\incl{\hookrightarrow}

\def\rad{\hbox{rad}}
\def\Spin{\hbox{Spin}}
\def\ov{\overline}

\def\wt{\widetilde}

\def\boite#1{\leavevmode\hbox{\vrule\vtop{\vbox{\hrule\kern1pt\hbox
{\kern1pt\strut#1\kern1pt}}
\kern1pt\hrule}\vrule}}

\def\n{\cap}  
\def\t{\mathop{\otimes}}  
\def\ds{\mathop{\oplus}}

\def\Tor{\mathop{\rm Tor}\nolimits}
\def\\{\backslash}
\def\back/{\backslash}

\def\phii#1{\f \hskip -2pt \raise -3pt \hbox{${}_#1$}}

%
\def\hfl#1#2{\smash{\mathop{\hbox to 12mm{\rightarrowfill}}
\limits^{\scriptstyle#1}_{\scriptstyle#2}}}
\def\vfl#1#2{\llap{$\scriptstyle#1$}\left\downarrow
\vbox to 6mm{}\right.\rlap{$\scriptstyle#2$}}
\def\hfle#1#2{\smash{\mathop{\hbox to 12mm{\leftarrowfill}}
\limits^{\scriptstyle#1}_{\scriptstyle#2}}}
\def\vfle#1#2{\llap{$\scriptstyle#1$}\left\uparrow
\vbox to 6mm{}\right.\rlap{$\scriptstyle#2$}}
\def\hfld#1#2{\smash{\mathop
{\buildrel {\hbox to 12mm{\rightarrowfill}} 
\over{{\hbox to 12mm{\leftarrowfill}}}}
\limits^{\scriptstyle#1}_{\scriptstyle#2}}} 
\def\diagram#1{\def\normalbaselines{\baselineskip=0pt
\lineskip=10pt\lineskiplimit=1pt} \matrix{#1}}
\def\vasur{\hbox{\raise -1pt\hbox{$\vdash$}}\kern -3pt\joinrel\mathrel{\hfl{}{}}}
\def\loong{\raise 2.2pt\vbox{\hrule width 2mm}\kern -3pt\joinrel\mathrel{\hfl{}{}}}
\def\longincl{ \lhook\joinrel\mathrel{\hfl{}{}} }
%
%
\cl{\dc Statement \ of \ the} \bs
\cl{\dc Alexandru \ Conjecture} \bs \bs
The purpose of this text is to add (at least conjecturally) some more items 
to the list 
of analogies between the category $\cal H$ of Harish-Chandra modules and the 
category $\cal O$ of Bernstein-Gelfand-Gelfand which has been established by 
Bernstein-Gelfand-Gelfand, Vogan, Beilinson, Ginzburg, Soergel and others. 
These 
analogies have been suggested by confronting some observations about 
$p$-integrable 
harmonic forms on real hyperbolic space with results of the people mentioned 
above 
about the category $\cal O$. 
To the reader more familiar with $L^p$ harmonic forms than with the category 
$\cal O$ 
my advice is to read in parallel Parts A and B of section~1 (Part B being a 
detailed 
example).  \bs \bs
\cl{\dz 1. The \ main \ statements} \bs \bs 
\cl{\uc Part A. The Weak Alexandru Conjecture} \bs
{\parindent=10mm
\item{(1.1)} {\bf Setting}.

\item{} $G$ is a connected semisimple Lie group $G$ with finite center,  

\item{} $K \subset G$ is a maximal compact subgroup, 

\item{} $\go g \supset \go k$ are the complexified Lie algebras of $G$ and 
$K$, 

\item{} $\go b \subset \go g$ is a Borel subalgebra, 

\item{} $\go h$ a Cartan subalgebra of $\go g$ contained in $\go b$. \par} \bs
For any pair $\go m \subset l$ of (complex) Lie algebras and any 
$\go l$-module $V$, say 
that $V$ is $\go m$-{\bf finite} if it is a sum of 
finite dimensional sub-$\go m$-modules, and that $V$ is an 
\hbox{$(\go l,m)$-{\bf module}} if 
it is $\go m$-finite and \hbox{$\go m$-semisimple.} The category $\cal O$ of 
\hbox{\bf BGG-modules} is the full subcategory of $\go g$-mod 
whose objects are the $\go b$-finite $(\go g,h)$-modules of finite length ; 
whereas  
the category $\cal H$ of {\bf Harish-Chandra modules} is the full subcategory of 
\hbox{$\go g$-mod} 
whose objects are those $(\go g,k)$-modules of finite length $V$ such that for 
any finite dimensional \hbox{$\go k$-invariant} subspace 
$F \subset V$ the action of $\go k$ on $F$ exponentiates to $K$. The 
categories $\cal O$ 
and $\cal H$ are $\C$-categories in the sense of Bass [B] page 57.  \bs
\ind For any $\C$-category $\cal C$ let
$$\cal I = I(C)$$
be the set of isomorphism classes of simple objects of $\cal C$ [assume it 
{\it is} a set], for each $i \in \cal I$ choose a representative
$$L_i \in i$$
and let
$$\l(i)$$
be the {\bf projective dimension} of $L_i$ [{\it i.e.} the supremum in 
$\Z \cup \{+ \infty\}$ of 
the set $\{n \in \Z \ | \ \Ext^n(L_i,-) \not= 0 \}$]. \bs 
%
%
\vbox{ {\parindent=10mm 
\item{(1.2)} {\bf Definition}. The {\bf $\cal C$-ordering} is the smallest partial 
ordering $\leq$ on $\cal I$ satisfying 
$$\left.\matrix{i, j \in {\cal I} \ \cr \cr
\l(j) = \l(i) + 1 < \infty \ \cr \cr
\Ext^1(L_j,L_i) \not= 0 \ \cr}
\right\} \ \Longrightarrow \ i \leq j.$$} \par} 
%
%
\vbox{{\parindent=10mm 
\item{(1.3)} {\bf Definition}. The {\bf subcategory generated by the subset} 
$\cal J$ {\bf of} $\cal I$ is the full sub-$\C$-category $\( \cal J \)_C$ of $\cal C$ 
characterized by the condition that an object $V$ of $\cal C$ belongs to 
$\( \cal J \)_C$ iff each simple subquotient of $V$ is isomorphic to $L_j$ for 
some $j \in \cal J$. \par}}   \bs
%
%
\vbox{{\parindent=10mm 
\item{(1.4)}  {\bf Definition}. If $\cal C$ is a $\C$-category and $\cal B$ 
a full sub-$\C$-category, say that $\cal B$  is {\bf Ext-full} in $\cal C$ if for all 
$V, W \in \cal B$ the natural morphism 
$$\Ext_{\cal B}^\.(V,W) \rightarrow \Ext_{\cal C}^\.(V,W)$$
is an isomorphism. \par}} \bs 
%
%
\ind [If all objects of $\cal B$ have finite length it suffices to check the above 
isomorphism for $V$ and $W$ simple (because of the long exact sequences and the 
five-lemma).] 
%
%
Recall that a subset $\cal J$ of $\cal I$ is an {\bf initial segment} iff 
$$\left.\matrix{i, j \in {\cal I} \ \cr \cr
j \in {\cal J} \ \cr \cr
i \leq j \ \cr}
\right\} \ \Longrightarrow \ V \in {\cal J}.$$ 
%
%
%
\vbox{ {\parindent=10mm 
\item{(1.5)} {\bf Definition}. In the above notation $\cal C$ is a {\bf Guichardet 
category} if the subcategory generated by any initial segment is Ext-full 
in $\cal C$. \par}} \bs
\ind Let Setting (1.1) be in force, denote by $I$ the annihilator of the trivial module 
in the center of $U(\go g)$ and for any sub-$\C$-category $\cal C$ of $\go g$-mod let 
$${\cal C}_\rho$$
be the full sub-$\C$-category of $\cal C$ whose objects are annihilated by some power 
of $I$.  \bs
%
%
\vbox{ {\parindent=10mm \item{(1.6)} {\bf Weak Alexandru Conjecture}. 
The category ${\cal H}_\rho$ is a Guichardet category. \par}} \bs
%
%
\ind Part (a) of the theorem below is due to Cline, Parshall and Scott, and 
\hbox{part (b)} to Fuser. \bs
%
%
\vbox{{\parindent=10mm \item{(1.7)} {\bf Theorem}. 

\itemitem{(a)} The category ${\cal O}_\rho$ is a Guichardet category,

\itemitem{(b)} the Weak Alexandru Conjecture holds for ${\cal H}_\rho$ with $G = SL(3,\R)$, 
$\Spin(n,1)$ or $SU(n,1)$. \par}} \bs
%
%
\ind The next item on the agenda is the Strong Alexandru 
Conjecture (SAC). Hoping to make it more digest I first give a set of statements 
of a somewhat geometric flavor which imply the SAC for \hbox{$\Spin(2n+1,1) \ $;} 
more precisely the setting of the SAC is dual to the one described here. I'll 
use horizontal lines to set off this part of the text [which is merely 
motivational]. \bs \hrule \bs
%
%
\cl{\uc Part B. A detailed example} \bs
Let $\Omega^p_k$ be the space of those $p$-forms on hyperbolic \hbox{$(2n+1)$-space} 
which are killed by $\Delta^k$ and $\Omega^p$ the union of the $\Omega^p_k$ ; put 
$${\cal I} := \{ p \in \Z \ | \ 0 \leq p \leq n \},$$
$$\Omega := \ds_{p \in \cal I} \Omega^p.$$
As a general notation if $V$ is a vector space, $W$ a vector subspace and 
$\f_1,..., \f_k$ endomorphisms of $V$, set 
$$W(\f_1,... ,\f_k) := W \n \left( \ \bigcap_{i=1}^k \Ker \f_i \right).$$
Let $d$ be the differential and $d^*$ be the codifferential, and put for 
$p \in \cal I$
$$ L_p := \Omega^p(d, d^*),$$
$$ M_p := \Omega^p(d^*),$$
$$\ov{M}_p := \left\{\matrix{
\Omega^p(\Delta, d^*) & \hbox{ if } & p < n ,  \cr \cr
L_n & \hbox{ if } & p = n .}\right.$$
Equip $\Omega^p_k$ with the $C^0$-topology and $\Omega^p$ with the inductive limit 
topology ; 
let $G$ be the group of orientation preserving hyperbolic isometries ; denote by 
$\cal C$ the category defined by the rule that an object of $\cal C$ is a topological 
$G$-module $V$ which is isomorphic to a close subspace of some 
$\Omega^{p_1} \ds \cdots \ds \Omega^{p_r}$, and 
a morphism in $\cal C$ is a \hbox{$G$-equivariant} continuous linear map. The following 
facts are known : \bs
{\parindent=11pt 
\item{$\.$} The category $\cal C$ is a $\C$-category.  \bs
\item{$\.$} The modules introduced above belong to it ; any simple object of $\cal C$ 
is isomorphic to $L_p$ for a unique $p$ in $\cal I$.  \bs
\item{$\.$} The $\cal C$-ordering on $\cal I$ is opposite to the 
natural ordering ; the projective dimension of $L_p$ is \hbox{$2n+1-p$.}  \bs
\item{$\.$} The category ${\cal H}_\rho$ is 
equivalent to the subcategory of $\cal C$ whose objects have finite length, or 
equivalently are annihilated by some power of $\Delta$.  \bs 
\item{$\.$} This subcategory --- which I abusively denote by 
${\cal H}_\rho$ for a short while --- is \hbox{Ext-full} in $\cal C$ and contains 
$\ov{M}_p \, $. \bs
\item{$\.$} The objects of $\cal C$ have injective hulls and therefore 
minimal injective resolutions. \par}  \bs
[In this parenthesis I give a reminder of what's meant here by minimal resolution and 
offer to the reader unfamiliar with homological algebra a cheap definition 
of \hbox{Ext-groups} in this context. Given $V$ in $\cal C$ there are elements 
$p_1,... ,p_k$ of $\cal I$ and an 
isomorphism from the socle of $V$ onto $\ds_{i=1}^k L_{p_i} \ ;$ this isomorphism 
extends to an embedding $\f : V \mono I^0 := \ds_{i=1}^k \Omega^{p_i}$ (this statement is 
sometimes called Frobenius reciprocity). Since the cokernel of $\f$ is 
again in $\cal C$ this process can be iterated, giving rise to an injective 
resolution $V \mono I^0 \arrow I^1 \arrow \cdots$, which is clearly minimal. Moreover 
for each pair of integers $(j,p)$ with $j \geq 0$ and $p \in \cal I$ there is 
a finite dimensional vector space $\Ext^j(L_p,V)$ acted on trivially by $G$ such that 
$I^j \simeq \ds_{p \in \cal I} \ \Ext^j(L_p,V) \t_\C \Omega^p$.] 
\bs
{\parindent=11pt 
\item{$\.$} The injective hull of $L_p$ is $\Omega^p$.  \bs
\item{$\.$} The vector spaces $\Ext^\.(V,W)$ are finite dimensional 
for $V, W \in \cal C$. \par} \bs
\ind Let me digress a tiny bit by stating the {\it Weak} Alexandru Conjecture in this 
setup. During this parenthetical comment $p$ shall be a fixed ``number" satisfying  
$1 \leq p \leq \infty$ [{\it not} an element of $\cal I$]. Let $K$ be the stabilizer 
in $G$ of a given point of hyperbolic space, say that an object $V$ of ${\cal H}_\rho$ is 
$p$-{\it integrable} if its \hbox{$K$-finite} vectors are, and let ${\cal L}^p$ be the full 
subcategory of $p$-integrable objects of ${\cal H}_\rho$ ; then ${\cal L}^p$ is Ext-full in 
${\cal H}_\rho$ [more precisely a subcategory of ${\cal H}_\rho$ is generated by an 
initial segment iff it is of the form ${\cal L}^p$ --- the number $p$ being of course 
in general nonunique.] \bs 
\ind Going back to the SAC, for $q \in \cal I$ 
consider the filtration $F_0^q := 0 \subset F_1^q := M_q \subset F_2^q := \Omega^q$ ; 
then $d^*$ induces an isomorphism $F_2^q / F_1^q \simeq M_{q-1} \ $ [with the convention 
$M_{-1} = 0$] ; this filtration is analogous to the 
filtration of projective modules by Verma modules in the category $\cal O \ $  
and is encoded in Axiom (1.10) below. If $E^\.$ is a graded vector space let 
$E^\.(t) $ be its Poincaré series. Define the 
$\cal I$ by $\cal I$ matrix $a$ with entries in $\Z[t]$ by
$$a_{pq} = \Ext^\.(L_p , M_q)(t).$$
To compute this series note that the augmented complex 
$$\diagram{
M_q & \longincl & \Omega^q & \hfl{d^*}{} & \Omega^{q - 1} & \hfl{d^*}{} &
\cdots & \hfl{d^*}{} & \Omega^0 \cr} $$ 
is ``the" minimal injective resolution of $M_p \,$, whence 
$$a_{pq} = \left\{\matrix{
t^{q-p} & \hbox{ if } & p \leq q ,  \cr \cr
0 & \hbox{ if } & q < p .}\right.$$
Letting $[V]$ be the class of $V$ in the Grothendieck group $[\cal C]$ of 
$\cal C$ we have 
$$\left[ \, \ov{M}_p \right] = [L_p] + [L_{p+1}]$$
[with $L_{n+1} = 0$] and 
$$[L_p] = \ \sum_{q=p}^n \ (-1)^{q-p} \ \left[ \, \ov{M}_q \right] = \ 
\sum_{q \in \cal I} \ a_{pq}(-1) \ \left[ \, \ov{M}_q \right]  \ ;$$
this corresponds to the Delorme formula in the category 
$\cal O$ and is encapsulated in Axiom (1.14). The inverse $a^{-1} = (a^{-1}_{pq})$ 
of $a$ being given by
$$a^{-1}_{pq} = \left\{\matrix{
(-t)^{q-p} & \hbox{if } p \leq q \leq p+1 ,  \cr \cr
0 & \hbox{otherwise}}\right.$$ 
the number of occurrences of $M_p$ in the filtration $F^q_\.$ is 
$a^{-1}_{pq}(-1)$, which is an analog of the BGG duality in the category 
$\cal O$, and gives rise to Axiom (1.15) below. 
Let $\ov{M}_{p,\hbox{\sevenrm soc}}^\.$ the graded object of 
$\cal C$ associated to the socle filtration of $\ov{M}_q$ and 
$\big[ \, \ov{M}_{q,\hbox{\sevenrm soc}}^\. \big](t)$ be its image in 
$\Z[t] \t_\Z [\cal C]$, then
$$\left[ \, \ov{M}_{p,\hbox{\sevenrm soc}}^\. \right](t) = [L_p] + t \ [L_{p+1}] = \ 
\sum_{q \in \cal I} \ a^{-1}_{pq} \ [L_q].$$
The corresponding theorem in the realm of the category $\cal O$ is due to 
Beilinson, Ginzburg and Soergel, and the above formula suggests Axiom 
(1.16) below. Finally let's compute $\Ext^\.(L_p,L_q)$. To simplify the notation,  
if $n < p \leq 2n+1$ identify $\Omega^p$ to $\Omega^{2n+1-p}$ via the star 
operator, and decree that $\Omega^p = 0$ for $p < 0$ or $p > 2n+1$. 
``The" minimal resolution of $L_q$ being
$$\diagram{
L_q & \incl & \Omega^q & \hfl{ \big( {\kern -4pt d \atop d^*} \big)}{} & 
\Omega^{q+1} \ds \Omega^{q-1} & 
\hfl{ \big( {\kern -7pt d \ 0 \atop \kern -1pt 0 \hskip 3pt d^*} \kern -2pt \big) }{} &
\Omega^{q+2} \ds \Omega^{q-2} & 
\hfl{ \big( {\kern -7pt d \ 0 \atop \kern -1pt 0 \hskip 3pt d^*} \kern -2pt \big)}{} & 
\cdots \cr}$$
we have
$$\Ext^\.(L_p,L_q)(t) = t^{|q-p|} + t^{2n+1-p-q}. \eqno{(*)}$$ 
Setting $L := \ds_{p \in \cal I} L_p \,$ and letting ${}^t \kern -2pt a$ be 
the transpose of $a$ and $\delta$ the diagonal matrix defined by
$$\delta_p = \left\{\matrix{
1+t & \hbox{if } p = n,  \cr \cr
1-t^2 & \hbox{if } p \not= n,} \right.$$
$(*)$ reads
$$\Ext^\.(L,L)(t) = a \ \delta \ {}^t \kern -2pt a.$$
The category $\cal O$ analog is the Beilinson-Ginzburg formula and the 
corresponding Axiom below is (1.18). \bs \hrule \bs
%
%
\cl{\uc Part C. The Strong Alexandru Conjecture} \bs
Here are some preliminaries to state the Strong Alexandru Conjecture. Set
$$Z := \C[[z_1,...,z_m]],$$
[where $z_1,..., z_m$ are indeterminates] and let $A$ be a $Z$-algebra which 
is 
finitely generated over $Z$. Then there is a semisimple subalgebra $A_0$ 
of $A$ satisfying \hbox{$A = A_0 \ds \rad(A)$.} Assume there is a finite set 
$F$ 
such that $A_0$ can --- and will --- be identified to the algebra $\C^F$ of 
functions on $F$. For each $i \in F$ define 
$e_i \in A_0$ by $e_i(j) = \delta_{ij}$ [Kronecker delta]. As a general 
notation put
   $$A\hbox{-fd} := \hbox{ the category of finite dimensional $A$-modules}.$$
%
%
I'll make free use of the facts that by a 
theorem of Casselman $A$-fd is Ext-full in \hbox{$A$-mod} and that by a 
result of 
BGG the categories ${\cal O}_\rho$ and ${\cal H}_\rho$ are equivalent to 
\hbox{$A$-fd} for some algebra $A$ as above --- the $\C$-algebra isomorphism 
class of $A$ being unique. 
It will be tacitly assumed that the category $\cal C$ of interest 
has been set to be ${\cal O}_\rho$ or ${\cal H}_\rho \,$, and that an 
algebra $A$ as above and an equivalence \hbox{${\cal C} \sim A$-fd} have been 
chosen, 
providing in particular an identification \hbox{${\cal I} = F$ ;} the symbol 
$L_i$ denotes at the same time an object of $\cal C$ and ``the" 
corresponding object in $A$-fd ; more generally I'll allow myself to 
navigate rather freely between $\cal C$ and $A$-fd. Fix a family 
$M = (M_i)_{i \in \cal I}$ of \hbox{$A$-modules.} \bs
%
%
\vbox{ {\parindent=9mm 
\item{(1.8)} {\bf Definition}. An {\bf $M$-filtration} $F_\.$ of an $A$-module $V$ is a 
sequence 
$$0 = F_0 \subset F_1 \subset F_2 \subset \cdots \subset F_r = V$$
of submodules such that there exists an $r$-tuple $(i_1,...,i_r)$ of 
elements of $\cal I$ satisfying $F_p/F_{p-1} \simeq M_{i_p}$ for 
$1 \leq p \leq r$. An $A$-module is $M$-{\bf filtrable} if it admits an 
\hbox{$M$-filtration.}   \par}}    \bs
%
%
\ind Let $\leq$ be the $\cal C$-ordering on $\cal I$ [see Definition (1.2)] ; form the 
small Verma module 
$$M^-_i(A) := Ae_i \Bigm/ \sum_{j \not\leq i} A \, e_j \, A \, e_i$$
and the large Verma module 
$$M^+_i(A) := Ae_i \Bigm/ \sum_{j > i} A \, e_j \, A \, e_i \ ;$$
letting $M_i$ be as in Definition (1.8) set
$$\ov{M}_i := M_i \bigm/ \rad(\End_AM_i) \, M_i \, ,$$
$$E^\.(t) := \hbox{Poincaré series of the graded vector space $E^\.$,}$$
$$[V] := \hbox{class of $V$ in the Grothendieck group.}$$
Let $a$ be the $\cal I$ by $\cal I$ matrix with entries in $\Z[[t]]$ defined by
$$a_{ij}(t) := \Ext_A^\.(M_i,L_j)(t)$$
and if $a$ happens to be invertible let $a^{-1}_{ij}$ be the $(i,j)$ entry of 
$a^{-1}$.
%
%
Let $V$ be in \hbox{$A$-fd ;} as a general notation set 
$$V_{\hbox{\sevenrm rad}}^\. := \ds \  
\left(\rad(A)^i \, V \bigm/ \rad(A)^{i+1} \, V \right) \ ;$$

suppose $V \not= 0$ ; let $S_{-1} := 0 \subset S_1 = $ soc $V \subset \cdots 
\subset S_p = V$ be the socle filtration of $V$, with $S_{p-1} \not= V$ ; 
say that {\bf the radical and socle filtrations coincide} if 
$$\rad(A)^{i} \hskip 1mm V = S_{p-i} \ \ \forall \ \ i.$$
%
Let ${}^t b$ denote the 
transpose of any matrix $b$ ; and consider the following conditions \bs
%
%
{\parindent=11mm 
\item{(1.9)} for each $i \in \cal I$ we have 
$M^+_i(A) = M^-_i(A) = M_i \ $;  \bs
\item{(1.10)} for each $i \in \cal I$ the module $Ae_i$ is 
\hbox{$M$-filtrable} [see Definition (1.8)] ;  \bs
\item{(1.11)} for each $i$ in $\cal I$ the module $M_i$ is flat over 
$\End_A M_i \ ;$ \bs
\item{(1.12)} for any $i \in \cal I$ there is a $j \leq i$ such that 
$$\l(i) = \sup \ \{ n \in \Z \ | \ \Ext^n(L_i \, , L_j) \not= 0 \} \ ;$$
%
%
\item{(1.13)} there are polynomials $p_{ij}$ such that 
$$a_{ij}(t) = t^{\l(j) - \l(i)} \ p_{ij}(t^{-2}),$$
$$p_{ij} \not= 0 \iff i \leq j \iff p_{ij}(0) = 1,$$
$$p_{ii} = 1,$$
$$\deg P_{x,y} \ < \ {\l(y)-\l(x) \over 2} \ \mbox{ if } \ 
x < y \ ;$$
in particular $a$ is invertible ; \bs
\item{(1.14)} $[L_j] = \sum_i \ a_{ij}(-1) \ \left[ \, \ov{M}_i \right]$ ; \bs
\item{(1.15)} things can be arranged so that the number of 
occurrences of $M_j$ in the \hbox{$M$-filtration} of $Ae_i$ in (1.10) is 
$a^{-1}_{ij}(-1)$  ;  \bs
\item{(1.16)} the radical and socle filtrations of $\ov{M}_j$ coincide 
and we have 
$$\left(e_i \hskip 2pt \ov{M}_{j,\hbox{\sevenrm rad}}^{ \ \.} \right)(t) = 
a^{-1}_{ij}(-t) \ ; $$
\item{(1.17)} $\End_A\big(\ov{M}_i\big) = \C$ \ ; \bs
\item{(1.18)} there is a diagonal $\cal I$ by $\cal I$ matrix $d$ such that 
$$\Ext^\._A(A_0,A_0)(t) = {}^ta \ d \ a.$$ \par}
%
%
\ind The proposition below is essentially due to Cline, Parshall and Scott ; 
an elementary proof is given in section~1 [recall that $A$-fd is the category of 
finite dimensional $A$-modules]. \bs 
%
%
\vbox{ {\parindent=11mm 
\item{(1.19)} {\bf Proposition}. In the above setting if Conditions (1.9)---(1.12) 
are satisfied then $A$-fd is a Guichardet category. \par}} \bs  
%
%
\vbox{{\parindent=11mm 
\item{(1.20)} {\bf Definitions}. If $A$ satisfies 
Conditions (1.9)---(1.18) above, then $A$ is a {\bf BGG algebra}. A 
{\bf beegeegee} is a category which is equivalent to $A$-fd for some BBG algebra 
$A$. \par}}  \bs
%
%
\ind The theorem below is due to BGG (see [BGG]), Beilinson and Ginzburg (see [BGS]) and 
Cline-Parshall-Scott (see statements (3.3.c), (3.5.a) and (3.9.a) in [CPS]). \bs
%
%
\vbox{(1.21) {\bf Theorem}. The category ${\cal O}_\rho$ is a beegeegee.}  \bs
%
%
\vbox{(1.22) {\bf Strong Alexandru Conjecture}. The category ${\cal H}_\rho$ is a 
beegeegee.}  \bs
%
%
\vbox{{\parindent=11mm 
\item{(1.23)} {\bf Theorem} (Fuser). The above conjecture holds for $\Spin(n,1), SU(n,1)$ 
and $SL(3,\R)$. \par}}  \bs
%
%
\ind The drawback [at least one of them] of all this stuff is that it's 
almost never computable ! Here is a statement which, although as conjectural 
as the previous ones, can be submitted to numerical tests. It consists in a 
computable 
variant of [a particular case of] Condition (1.18), that is in a formula which 
would express, when $\go g$ and $\go k$ have the same rank, 
the Poincaré series $\Ext^\._A(A_0,A_0)(t)$ in terms of computable things. 
In the Langlands classification $L_i$ occurs as the unique 
simple quotient of a module induced from some parabolic subgroup $P_i$ ; let 
\hbox{${\go p}_i = {\go m}_i \ds {\go a}_i \ds {\go n}_i$} be 
a Langlands decomposition of $\C \t_\R \hbox{Lie}(P_i)$ ; put
$$\wt d_i := \left(1-t^2 \right)^{\dim {\go a}_i} \ ;$$
%
%
let $\wt \l(i)$ be the dimension of 
the $K_\C$-orbit attached to $i$ and $(\wt p_{ij})$ the family of 
Kazhdan-Lusztig-Vogan polynomials, that is 
the one denoted $(P_{\gamma, \delta})$ in [V3], section~6 ; and set 
$$\wt a_{ij}(t) = t^{\wt \l(j) - \wt \l(i)} \ \wt p_{ij}(t^{-2}) \, .$$
%
%
(1.24) {\bf Conjecture}. If $\go g$ and $\go k$ have the same rank then
$$\Ext^\._A(A_0,A_0)(t) = {}^t \wt a \ \wt d \ \wt a.$$
%
%
\ind For a given group $G$ one can run the following numerical test. \bs

\vbox{\parindent=6mm \item{(a)} Compute ${}^t \wt a \ \wt d \ \wt a$ ;

\item{(b)} look if this is compatible with what's known of 
$\Ext^\._A(A_0,A_0)(t)$ [in 
particular with the $({\go g},K)$-cohomology, as computed by Vogan's 
\hbox{$U_\alpha$-algorithm}] ;

\item{(c)} pretend Conjecture (1.24) holds and use it to compute the 
${\cal H}_\rho$-ordering ; 

\item{(d)} check if there are matrices $a$ and $d$ such that $a$ is upper 
triangular with ones on the diagonal, $d$ is diagonal and we have 
${}^t a \ d \ a = {}^t \wt a \ \wt d \ \wt a$. \par} \bs
%
%
\ind If the test is successful no conclusion can be drawn ; if it fails then 
at least one of the involved conjectures is wrong. --- In the real rank one case 
evidence suggests that the classical objects ``with tildes" coincide with the 
nonclassical ones [``without tildes"] ; in the case of $PSp(2,\R)$ using results of 
Vogan [V], \hbox{pp 251-255,} one sees such is {\bf not} the case --- but the test is still 
successful. I think that in the case where the complex ranks of $\go g$ and $\go k$ 
are different there is a similar formula with $\wt a$ as above and $\wt d$ a certain 
polynomial valued diagonal matrix. \vfill\eject
%
%
%
%
\cl{\dz 1. \ Proof \ of \ Proposition \ (1.19) } \bs \bs
In the setting of Proposition (1.19) assume 
there are at least two elements in $\cal I$ [otherwise there is 
nothing to prove], and consider the following setting 
\vbox{$$i \hbox{ is a maximal element of } {\cal I},$$ 
$$e := e_i \, ,$$
$$I := AeA,$$
$$B := A/I,$$
$${\cal J} := {\cal I} \ \\ \ \{i\}.$$}
Proposition (1.19) follows from Lemma (2.1) below, which will be proved at 
the end of the section. \bs
%
%
\vbox{{\parindent=9mm 
\item{(2.1)} {\bf Lemma}.  \smallskip
\itemitem{(a)} The category $B$-mod is Ext-full in $A$-mod, \smallskip
\itemitem{(b)} The $\cal C$-ordering coincides on $\cal J$ with the 
$\( \cal J \)_C$-ordering [see (1.2) and (1.3)],  \smallskip
\itemitem{(c)} $B$ satisfies Conditions (1.9)---(1.12). \par}} \bs
%
%
\vbox{(2.2) {\bf Lemma}. If $j \not= i$ then $M_j \in B$-fd.} \bs
%
%
{\bf Proof of Lemma (2.2)}. The statement is an immediate consequence of the 
following observation. For any pair $({\cal K},k)$ with 
$j \in \cal K \subset I$ put
%
%
$$P({\cal K},k) :=  Ae_k \Bigm/ \sum_{\l \not\in {\cal K}} A \, e_\l \, A \, e_k \, .$$  
Then $P({\cal K},k)$ is a projective cover of $L_k$ in 
$\({\cal K} \)_{A \hbox{\sevenrm -fd.}}$ If 
\hbox{$k \in \cal L \subset K \subset I$} and $P({\cal K},k)$ is in 
$\( {\cal L} \)_{A \hbox{\sevenrm -fd}}$ then $P({\cal K},k) = P({\cal L},\l)$. 
{\bf QED}  \bs
%
\ind The following lemma is obvious.   \bs
%
%
\vbox{{\parindent=9mm 
\item{(2.3)} {\bf Lemma}. The three conditions  
(a) $j$ is maximal, \smallskip

\hskip 49mm (b) $M_j = Ae_j \,$, \smallskip

\hskip 49mm (c) $M_j$ is projective  \smallskip

are equivalent. In particular $M_i = Ae_i$ is projective. {\bf QED} \par}} \bs
%
%
\ind Let $\cal E$ be the class of those $A$-modules which are isomorphic to a 
direct sum of finitely many copies of $Ae_i \,$.  \bs 
%
%
\vbox{ {\parindent=9mm \item{(2.4)} {\bf Lemma}. Let $V$ be an $M$-filtrable 
$A$-module. Then there is an \hbox{$M$-filtration} $F_\.$ of $V$ and a nonnegative 
integer $n$ such that \hbox{$F_n \in \cal E$} and $e(V/F_n) = 0$. \par}}  \bs 
%
%
{\bf Proof}. The statement results from the fact easy to check that 
\hbox{$eM_j = 0$} for all $j \not= i$ and from the observation that 
each occurrence of a projective into a given module as a subquotient is in
fact an occurrence as a submodule. {\bf QED} \bs
%
%
\ind In particular the above lemma provides for each $j \in \cal I$ an 
{$M$-filtration} $F^j_\.$ of $Ae_j$ by left subideals and a nonnegative integer 
$n_j$ such that \hbox{$F^j_{n_j} \in \cal E$} and \hbox{$e(Ae_j /F^j_{n_j}) = 0$.} Put 
%
%
$$J := \ds_j F^j_{n_j} \ ;$$
then there is a positive integer $n$ such that $J$ is isomorphic to the direct sum 
of $n$ copies of $Ae_i$ ; write this direct sum in the nonsensical form 
$\ds_{j=1}^n Ae_i$ ; choose an isomorphism 
$$\f : \ds_{j=1}^n Ae_i \iso J \ ;$$
%
%
and set 
$$A_i := \End_A Ae, $$
$$C := e A e = (A_i)^{op}.$$

We may --- and will --- view any left module over $A_i$ as a right module over 
$C$. Put 
%
%
$$\t := \t_C . $$
%
%
\vbox{ {\parindent=9mm \item{(2.5)} {\bf Lemma}. 
\itemitem{(a)} For \hbox{$j \not= i$} the multiplication map 
$Ae \t eAe_j \arrow Ae_j$ is an isomorphism, 
\itemitem{(b)} the natural map $Ae \t eA \arrow A$ is one-to-one,
\itemitem{(c)} we have $J = I$ and $F^j_{n_j} = Ie_j \,$. 
In particular $Be_j := Ae_j/Ie_j = Ae_j/F^j_{n_j}$ is $M$-filtrable. \par}}    \bs
%
%
{\bf Proof}. To prove (a) set \hbox{$f := e - e_j$} and note that the canonical 
isomorphism $Ae \t eAe \arrow Ae$ is the direct sum of the 
multiplication maps $Ae \t eAe_j \arrow Ae_j$ and 
$Ae \t eAf \arrow Af$. To check (b) consider the diagram
%
%
\def\dsp{\displaystyle}   
$$\diagram{ 
    Ae \t eA         & \hfl{\alpha_1}{} &   A  
\cr
    \vfl{\alpha_2}{} &                  &       \vfle{}{\f} 
\cr
    Ae \t eJ         &                  & \dsp{\ds_{j=1}^n Ae_i}
\cr
    \vfl{\alpha_3}{} &                  &       \vfle{}{\alpha_5}
\cr
    \dsp{Ae \t \ \ds_{j=1}^n eAe_i} 
                     & \hfl{}{\alpha_4} & 
                                    \dsp{\ds_{j=1}^n Ae \t Ae_i} \ ,
\cr}$$ 
%
%
where $\alpha_1,..., \alpha_5$ are defined as follows  \bs

{\parindent=20mm  $\alpha_1$ is the multiplication map, 

$\alpha_2$ is the identity [indeed $e(A/J) = 0 \then eA = eJ$], 

$\alpha_3$ is $1 \t \f^{-1}$,

$\alpha_4$ is the canonical map,

$\alpha_5$ is the multiplication map. \par} \bs
The maps $\alpha_2, \alpha_3$ and $\alpha_4$ are clearly isomorphisms and 
by (a) so is $\alpha_5 \,$, whereas $\f$ is one-to-one. Since the diagram 
commutes $\alpha_1$ is one-to-one and (b) is verified. To prove (c) 
observe $I = \Im \alpha_1 = \Im \f = J$.  {\bf QED} \bs
%
%
\vbox{ { \parindent=9mm \item{(2.6)} {\bf Lemma}. Let $A$ be an algebra and $e$ 
an element of $A$ satisfying $e^2 = e \not= 1$. Put $I := AeA$, $B := A/I$ and 
$C := eAe$, let $V$ and $W$ be $B$-modules,
and consider the following conditions  \medskip

\itemitem{(a)} the natural map $Ae \t_C eA \arrow I \subset A$ is 
one-to-one and $Ae$ is right  $C$-flat,  \smallskip
\itemitem{(b)} $I$ is right $A$-flat,  \smallskip
\itemitem{(c)} $\Tor_{q}^{A}(B,V) = 0$ for all $q > 0$,  \medskip
\itemitem{(d)} $\Ext^n_B(V,W) \simeq \Ext^n_A(V,W)$ for all $n$ (natural isomorphism).  
\medskip
Then (a) $\then$ (b) $\then$ (c) $\then$ (d). \par} } \bs
%
%
{\bf Proof}. The implication (a) $\then$ (b) is clear. The implication (b) $\then$ (c) 
follows from the long (?) exact sequence obtained by applying 
``$- \t_{A} V$" to the short exact sequence $I \mono A \epi B$. 
The implication  (c) $\then$ (d) follows from Proposition VI.4.1.3 of 
Cartan-Eilenberg [CE]. {\bf QED}  \bs
{\bf Proof of Lemma (2.1)}. Part (a) [the Ext-fullness of $B$-mod 
in $A$-mod] follows from Lemma (2.6). In view of Condition (1.12) [about projective 
dimensions] part (b) [the 
coincidence of the orderings] results from (a). Let me prove part (c), claiming  
that $B$ satisfies Conditions (1.9)---(1.12). Condition (1.9) [the 
coincidence of the Verma modules] and Condition (1.11) [the 
$(\End_B M_i)$-flatness of $M_i \,$] are consequences of \hbox{part (b)} and Lemma (2.2) 
\hbox{[$M_j \in B$-mod] ;} Condition (1.10) [the \hbox{$M$-filtrability]} follows from 
\hbox{(2.5.c) ;} Condition (1.12) results from (a). 
{\bf QED}\bs \cl{* \ * \ *} \bs 
%
%
%
%
{\parindent=15mm \lineskiplimit=10pt \lineskip=10pt
\item{[B]} Bass H., {\bf Algebraic K-theory}, Benjamin, New York 1968. \bs
\item{[BGG]} Bernstein I.N., Gelfand I.M., Gelfand S.I., Category of {\go g}-modules, 
{\it Funct. Anal. Appl.} {\bf 10}, 87-92 (1976).  \bs
\item{[BGS]} Beilinson A., Ginzburg V., Soergel W., Koszul duality patterns in 
representation theory, {\it J. Am. Math. Soc.} {\bf 9}, No.2, 473-527 (1996). \bs
\item{[CE]} Cartan H., Eilenberg S., {\bf Homological algebra}, Princeton University 
Press, 1956. \bs
\item{[CPS]} Cline E., Parshall B., Scott L., Finite dimensional algebras and 
highest weight categories, {\it J. Reine Angew. Math.} {\bf 391}, 85-99 (1988). \bs
\item{[D]} Delorme P., Extensions dans la cat\'egorie $\cal O$ de 
Bernstein-Gelfand-Gelfand. Applications, manuscript, October 1978. \bs
\item{[V]} Vogan D., The Kazhdan-Lusztig Conjecture for real reductive groups, in 
{\bf Representation theory of reductive groups}, Proceedings of the University of 
Utah conference 1982, Trombi Peter C. (Ed.), Birkhäuser, 1983, Progress in 
mathematics Number 40. \bs
\item{[V3]} Vogan, D., Irreducible characters of semisimple Lie groups. III : Proof 
of Kazhdan-Lusztig conjecture in the integral case, {\it Invent. Math.} {\bf 71} 
(1983) 381-417. \par} \bigskip

This text and others are available at http://www.iecn.u-nancy.fr/$\sim$gaillard
\bye